\documentclass[12pt]{amsart}

\usepackage{amsmath}
\usepackage{amsthm}
\usepackage{hyperref}

\numberwithin{equation}{section}

\newtheorem{prop}{Proposition}
\newtheorem{lemma}[prop]{Lemma}
\newtheorem{thm}[prop]{Theorem}

\numberwithin{prop}{section}

\theoremstyle{definition}
\newtheorem{defn}[prop]{Definition}
\newtheorem{ex}[prop]{Example}

\newtheorem{question}[prop]{Question}

\newcommand{\dt}{\frac{\partial}{\partial t}}
\newcommand{\brs}[1]{\left| #1 \right|}
\newcommand{\gG}{\Gamma}
\newcommand{\gD}{\Delta}
\newcommand{\gd}{\delta}
\newcommand{\gs}{\sigma}

\newcommand{\gl}{\lambda}
\newcommand{\gw}{\omega}
\newcommand{\ga}{\alpha}

\renewcommand{\ge}{\epsilon}
\newcommand{\N}{\nabla}

\renewcommand{\bar}[1]{\overline{#1}}
\newcommand{\del}{\partial}
\newcommand{\delb}{\bar{\partial}}
\newcommand{\bi}{\bar{i}}
\newcommand{\bj}{\bar{j}}
\newcommand{\bk}{\bar{k}}
\newcommand{\bl}{\bar{l}}
\newcommand{\bm}{\bar{m}}
\newcommand{\bp}{\bar{p}}
\newcommand{\bn}{\bar{n}}
\newcommand{\bq}{\bar{q}}
\newcommand{\br}{\bar{r}}
\newcommand{\bs}{\bar{s}}

\newcommand{\til}[1]{\widetilde{#1}}

\newcommand{\seven}{\mbox{VII}}

\newcommand{\LL}{\mathcal L}
\newcommand{\cf}{\frac{\sqrt{-1}}{2} \del \delb \log \det g}

\newcommand{\static}{static\ }

\DeclareMathOperator{\Sym}{Sym}

\DeclareMathOperator{\tr}{tr}
\DeclareMathOperator{\divg}{div}
\DeclareMathOperator{\Vol}{Vol}

\begin{document}

\title[A parabolic flow of pluriclosed metrics]{A parabolic flow of pluriclosed
metrics}

\author{Jeffrey Streets}
\address{Fine Hall\\
         Princeton University\\
         Princeton, NJ 08544}
\email{\href{mailto:jstreets@math.princeton.edu}{jstreets@math.princeton.edu}}

\author{Gang Tian}
\address{Fine Hall\\
	 Princeton University\\
	 Princeton, NJ 08544}
\email{\href{mailto:tian@math.princeton.edu}{tian@math.princeton.edu}}

\thanks{The first author was supported by the National Science Foundation via
DMS-0703660}
\thanks{The second author was partly supported by the National Science
Foundation via
DMS-0703985 and DMS-0804095}

\begin{abstract} We define a parabolic flow of pluriclosed
metrics.  This flow is of the same family introduced by the authors in
\cite{ST}.  We study the relationship of the existence of the flow and
associated static metrics topological information on the
underlying complex manifold.  Solutions to the \static equation are
automatically Hermitian-symplectic, a condition we define herein.  These static
metrics are 
classified on K3 surfaces, complex toroidal surfaces, nonminimal Hopf surfaces,
surfaces of general type, and Class $\seven^+$ surfaces.  To finish we discuss
how the flow may potentially be used to study the topology of Class $\seven^+$
surfaces.
\end{abstract}

\date{\today}



\maketitle

\section{Introduction}

In \cite{ST} the authors introduced a class of elliptic equations for
Hermitian metrics with associated parabolic flows.  A particular equation of
this type was singled out as being the unique such elliptic equation arising as
the Euler-Lagrange equation of a functional.  The purpose of this article is to
identify another equation in this same class related to pluriclosed metrics.
\begin{defn} Let $(M^{2n}, g, J)$ be a complex manifold with Hermitian metric
$g$ and K\"ahler form
$\omega = g(J \cdot, \cdot)$.  We say that $\omega$ is \emph{pluriclosed} if
\begin{gather} \label{pluriclosed}
\del \delb \omega = 0.
\end{gather}
\end{defn}

A result of Gauduchon \cite{Gauduchon} says that every conformal class of a
Hermitian metric on a complex $2n$-manifold has
a unique element satisfying
\begin{align*}
\del \delb \omega^{n - 1} = 0.
\end{align*}
We say that a metric satisfying this condition has \emph{null-eccentricity}.  In
the
case $n = 2$ the null-eccentricity condition agrees with the pluriclosed
condition.  Therefore every complex surface admits pluriclosed metrics.  The
condition
is considerably more restrictive in higher dimensions, though it is still in a
sense a determined equation for a complex $6$-manifold.  Some of the results in
this paper will apply to pluriclosed metrics in any dimension, but our main
focus will be on complex surfaces.

Before writing the definition of the flow let us briefly recall some definitions
related to the Chern connection.  In particular, let $\Omega$ denote the
curvature of the Chern connection $\N$ associated to $g$ and let
\begin{gather*}
S_{k \bl} = g^{i \bj} \Omega_{i \bj k \bl}.
\end{gather*}
Further, let
\begin{gather*}
T_{i j \bk} = \del_i g_{j \bk} - \del_j g_{i \bk}
\end{gather*}
denote the torsion of $\N$ and define
\begin{gather} \label{Qdefs}
\begin{split}
Q^1_{i \bj} =&\ g^{k \bl} g^{m \bn} T_{i k \bn} T_{\bj \bl m}\\
Q^2_{i \bj} =&\ g^{k \bl} g^{m \bn} T_{\bl \bn i} T_{k m \bj}.
\end{split}
\end{gather}
The tensor $Q^2$ will be used in certain evolution equations we derive later. 
In this
paper we will study the evolution equation
\begin{gather} \label{gflow}
\begin{split}
\dt g =&\ - S + Q^1\\
g(0) =&\ g_0.
\end{split}
\end{gather}
We note here that for certain applications a volume normalized version of
(\ref{gflow}) will be useful.  In particular one can add a scalar term to the
evolution equation to fix the volume of the manifold.  Setting $s = \tr_g S$ and
noting that $\tr_g Q^1 = \brs{T}^2$, this equation takes the form
\begin{gather} \label{vnflow}
\begin{split}
\dt g =&\ - S + Q^1 + \frac{1}{n} \left( \int_M s - \brs{T}^2 \right) g\\
g(0) =&\ g_0.
\end{split}
\end{gather}

\noindent These equations are of the form studied in \cite{ST}, and so the
regularity theory
developed in that paper applies immediately to (\ref{gflow}).
We will see in section 3 that this flow preserves the pluriclosed condition. 
Furthermore, as observed in \cite{ST},
if the initial condition is K\"ahler the resulting family of metrics is a
solution to K\"ahler-Ricci flow.  Moreover, in section 2 we will write this
equation using Hodge-type operators.  In particular a solution to (\ref{gflow})
with pluriclosed initial condition is equivalent to a solution of
\begin{gather} \label{flow}
\begin{split}
\dt \omega =&\ \del \del^* \omega + \delb \delb^* \omega
+ \frac{\sqrt{-1}}{2} \del \delb \log \det g\\
\omega(0) =&\ \omega_0.
\end{split}
\end{gather}
We prove the following basic regularity theorem in section 3.
\begin{thm} \label{regularity} Let $(M^{2n}, g, J)$ be a compact complex
manifold with pluriclosed
metric $g$.  There exists a constant $c(n)$ depending only $n$ such that there
exists a unique solution $g(t)$ to (\ref{gflow})
for
\begin{align*}
t \in \left[0, \frac{c(n)}{\max \{ \brs{\Omega}_{C^0(g_0)}, \brs{\N
T}_{C^0(g_0)}, \brs{T}^2_{C^0(g_0)}
\}} \right].
\end{align*}
Moreover, there exist constants $C_m$ depending only on $m$ such that the
estimates
\begin{align*}
\brs{\N^m \Omega}_{C^0(g_t)}, \brs{\N^{m+1} T}_{C^0(g_t)} \leq \frac{C_m \max \{
\brs{\Omega}_{C^0(g_0)}, \brs{\N T}_{C^0(g_0)},
\brs{T}^2_{C^0(g_0)} \}}{t^{m / 2}}
\end{align*}
hold for all $t$ in the above interval.  Moreover, if $g(0)$ is pluriclosed the
metric $g(t)$ is pluriclosed for all $t$ and is a solution to (\ref{flow}).  If
furthermore $g(0)$ is K\"ahler, then $g(t)$ is K\"ahler for all time and $g(t)$
solves K\"ahler-Ricci flow.
\end{thm}
\noindent A further consequence of this regularity is the following basic
long-time existence obstruction.
\begin{thm} Let $(M^{2n}, g(t), J)$ be a solution to
(\ref{gflow}).  Let $\tau$
denote the maximal existence time of the flow.  If $\tau < \infty$, then
\begin{align*}
\limsup_{t \to \tau} \max \left\{ \brs{\Omega}_{C^0(g_t)}, \brs{\N
T}_{C^0(g_t)}, \brs{T}^2_{C^0(g_t)} \right\} = \infty.
\end{align*}
\end{thm}

\noindent As it turns out, this regularity result can be significantly improved
in the case where $n = 2$.  In particular carefully analyzing the evolution of
the torsion and its covariant derivative on a complex surface allows us to
conclude that a bound on the Chern curvature suffices to show existence of the
flow.
\begin{thm} Let $(M^4, g(t), J)$ be a solution to (\ref{gflow}).  Let $\tau$
denote the maximal existence time of the flow.  If $\tau < \infty$, then
\begin{align*}
\limsup_{t \to \tau} \brs{\Omega}_{C^0(g_t)} = \infty.
\end{align*}
\end{thm}

Next we classify static solutions of (\ref{flow})
on certain complex surfaces.  In particular, given $(M^{2n}, g, J)$ a complex
manifold with pluriclosed metric, we say that $g$ is \emph{static} if
\begin{gather} \label{staticeqn}
- \del \del^* \omega - \delb \delb^* \omega
- \frac{\sqrt{-1}}{2} \del \delb \log \det g = \gl \omega.
\end{gather}
We note that if $g$ is K\"ahler and \static then it is K\"ahler-Einstein.  The
reason
for the sign on the left hand side is to make the sign of $\lambda$ agree with
the usual sign for KE-metrics.  On a complex surface the existence and signs of
\static metrics are closely related to the algebraic and topological structure
of
the manifold.  As we will see in section 4, \static metrics with nonzero
constant automatically imply the existence of a structure we call
Hermitian-symplectic.  Let us define this condition and expound on it.
\begin{defn} Consider $(M^{2n}, J)$ a complex manifold.  A
\emph{Hermitian-symplectic} \emph{structure} on $M$ is a real two-form $\omega$
such
that
$d \omega = 0$, and $\omega^{(1,1)} > 0$, i.e. the projection of $\omega$ onto
$(1,1)$-tensors determined by $J$ is positive definite.  We say that a complex
manifold is \emph{Hermitian-symplectic} if
it admits a Hermitian-symplectic structure.
\end{defn}
This shows that \static metrics do indeed carry a lot of structure with them. 
It is known that the space of symplectic manifolds is strictly larger than the
space of K\"ahler manifolds.  However, we do not know of an example or proof to
see if the space of Hermitian-symplectic manifolds is strictly larger
than the space of K\"ahler manifolds.  Indeed, one can make the following
observation.
\begin{prop} A complex surface is Hermitian-symplectic if and only if it is
K\"ahler.
\begin{proof} A K\"ahler structure is automatically Hermitian-symplectic, which
shows one direction.  Suppose $(M^{4}, J)$ is a complex surface and $\omega$ is
a Hermitian-symplectic form on $M$.  Suppose for a contradiction that $M$ is not
a K\"ahler manifold.  By the signature theorem the intersection form of a
non-K\"ahler surface is negative definite.  However, $\omega$ represents an
element of $H^2(M)$ with positive self-intersection.  This is a contradiction.
\end{proof}
\end{prop}
\noindent Therefore, any example of a complex manifold admitting a
Hermitian-symplectic structure but no K\"ahler structure must be of dimension
higher than $2$.  We state this question formally for emphasis.

\begin{question} Do there exist complex manifolds $(M^{2n}, J), n \geq 3$, such
that $M$ carries a Hermitian-symplectic structure $\omega$ but no K\"ahler
structure?
\end{question}

Let us now discuss the main results on static metrics.  As consequences of the
structural results in section 4 we give various
classification results for static metrics.  In particular, we show that static
metrics on K3 surfaces, two-dimensional complex tori, and surfaces of general
type must automatically be K\"ahler-Einstein.  Finally we can completely
determine the question of existence of static
metrics on Class $\seven^+$ surfaces (i.e. $b_1(M) = 1$ with $b_2(M) > 0$).
\begin{thm} \label{class7} Let $(M^4, J)$ be a complex surface of Class
$\seven^+$.  Then $M$ does not admit a
\static
metric.
\end{thm}
\noindent This theorem has important consequences in applying equation
(\ref{gflow}) to the study of the topology of class $\seven^+$ surfaces.  We
will discuss this in section 6.

Here is an outline of the rest of the paper.  In section 2 we give some general
background calculations for Hermitian metrics.  In section 3 we record the basic
existence and regularity results for (\ref{flow}), which are mainly consequences
of the general regularity theory developed in \cite{ST}.  Also we record a
number of evolution equations for various integral quantities along solutions to
(\ref{flow}).  In section 4 we give the improved regularity results in the case
of a complex surface.  Next in section 5 we derive certain equations and
inequalities
satisfied by \static metrics on surfaces, and classify static metrics on certain
complex surfaces.  Finally in section 6 we recall a basic
structural theorem for Class $\seven^+$ surfaces and give the proof of Theorem
\ref{class7}.

\section{Background Calculations}
In this section we fix notation and provide some basic calculations for standard
objects related to Hermitian geometry.  Fix $(M^{2n}, g, J)$ a complex manifold
with Hermitian metric $g$.  Let
\begin{gather*}
\gw(u,v) = g(J u, v)
\end{gather*}
be the K\"ahler form of $g$.  In local complex coordinates we have
\begin{gather*}
\gw = \frac{\sqrt{-1}}{2} g_{i \bj} dz^i \wedge d \bar{z}^j.
\end{gather*}
Let
\begin{gather*}
\Lambda^k = \bigoplus_{p+q = k} \Lambda^{p,q}
\end{gather*}
denote the usual decomposition of complex differential two-forms into forms of
type $(p,q)$.  The exterior differential $d$ decomposes into the operators
$\del$ and $\delb$
\begin{align*}
\del &: \Lambda^{p,q} \to \Lambda^{p+1, q}\\
\delb &: \Lambda^{p,q} \to \Lambda^{p, q+1}.
\end{align*}
Also the operator $d^*_{g}$, the $L^2$ adjoint of $d$, decomposes into
$\del^*_{g}$ and $\delb^*_{g}$
\begin{align*}
\del^*_g :&\ \Lambda^{p+1}, q \to \Lambda^{p, q}\\
\delb^*_g :&\ \Lambda^{p, q+1} \to \Lambda^{p,q}.
\end{align*}

\begin{lemma} \label{operatorlemma0} Given $g$ a Hermitian metric we have in
complex coordinates
\begin{align}
\left(\del^*_{g} \gw \right)_{\bk} =&\ \frac{\sqrt{-1}}{2} g^{p \bq}
\left(\del_{\bq} g_{p
\bk} - \del_{\bk} g_{p \bq} \right),\\
\left( \delb_g^* \gw \right)_j =&\ \frac{\sqrt{-1}}{2} g^{p \bq} \left(\del_j
g_{p \bq} - \del_p
g_{j \bq} \right).
\end{align}
\begin{proof} We compute using integration by parts.  Given $\ga \in
\Lambda^{0,1}$ we have
\begin{align*}
\left( \del_{g}^* \gw, \ga \right) =&\ \left( \gw, \del \ga \right)\\
=&\ \int_M g^{\bk l} g^{\bi j} \left(\gw_{j \bk} \bar{ \del \ga_{i \bl}} \right)
\bar{g}\\
=&\ \frac{\sqrt{-1}}{2} \int_M g^{\bi l} \left( \bar{\ga_{\bl, i}} \right)
\bar{g}\\
=&\ - \frac{\sqrt{-1}}{2} \int_M \bar{\ga_{\bl}} \left[ \del_{\bi} \left( g^{\bi
l} \bar{g} \right) \right]\\
=&\ - \frac{\sqrt{-1}}{2} \int_M \bar{\ga_{\bl}} \left(\bar{g}\right) \left[ -
g^{\bi m} \del_{\bi}
g_{m \bn} g^{\bn l} + g^{\bi l} \frac{1}{\bar{g}} \del_{\bi} \bar{g} \right].
\end{align*}
This gives the first formula, and the second follows analogously.
\end{proof}
\end{lemma}

\begin{lemma} \label{operatorlemma1} Given $g$ a Hermitian metric we have in
complex coordinates
\begin{align}
\left( \del \del^*_{g} \gw \right)_{j \bk} =&\ \frac{\sqrt{-1}}{2} \left[ g^{p
\bq} \left( g_{p \bk, \bq j} - g_{p \bq, \bk
j} \right) - g^{p \bm}
g^{n \bq} g_{\bm n, j} \left( g_{p \bk, \bq} - g_{p \bq, \bk} \right) \right]
\end{align}
and
\begin{align*}
\left(\delb \delb^*_g \omega \right)_{j \bk} =&\ \frac{\sqrt{-1}}{2} \left[ g^{p
\bq} \left(g_{j \bq, p \bk} - g_{p \bq, j \bk} \right) - g^{p \bm} g^{n \bq}
g_{n \bm, \bk} \left( g_{j \bq, p} - g_{p \bq, j} \right) \right].
\end{align*}
\begin{proof} In general for $\ga \in \Lambda^{0, 1}$ we have
\begin{align*}
\left(\del \ga \right)_{j \bk} = \del_j \ga_{\bk}.
\end{align*}
Thus we compute using Lemma \ref{operatorlemma0}
\begin{align*}
\left(\del \del_{g}^* \gw \right)_{j \bk} =&\ \frac{\sqrt{-1}}{2} \del_j \left(
g^{p \bq}
\left(\del_{\bq} g_{p \bk} - \del_{\bk} g_{p \bq} \right) \right)\\
=&\ \frac{\sqrt{-1}}{2} \left[ g^{p \bq} \left( g_{p \bk, \bq j} - g_{p \bq, \bk
j} \right) - g^{p \bm}
g_{\bm n, j} g^{n \bq} \left( g_{p \bk, \bq} - g_{p \bq, \bk} \right) \right].
\end{align*}
The result follows.  The second follows analogously.
\end{proof}
\end{lemma}

\begin{lemma} \label{operatorlemma10} Given $g$ a Hermitian metric we have in
complex coordinates
\begin{gather*}
\left(\frac{\sqrt{-1}}{2} \del \delb \log \det g \right)_{j \bk} =
\frac{\sqrt{-1}}{2} \left(g^{p \bq} g_{p \bq, j \bk} - g^{p \br} g^{s \bq}
g_{\br s, j} 
g_{p \bq, \bk} \right)
\end{gather*}
\begin{proof}
We compute directly in coordinates
\begin{align*}
\left(\frac{\sqrt{-1}}{2} \del \delb \log \det g \right)_{j \bk} =&\
\frac{\sqrt{-1}}{2} \del_{j}
\left( g^{p \bq}
\del_{\bk} g_{p \bq} \right)\\
=&\ \frac{\sqrt{-1}}{2}  \left(g^{p \bq} \del_j \del_{\bk} g_{p \bq} - g^{p \br}
\del_j g_{\br s} g^{s \bq}
\del_{\bk} g_{p \bq} \right).
\end{align*}
\end{proof}
\end{lemma}

In addition to the Hodge operators we will use two of the Ricci-type curvatures
of the Chern curvature.  In particular, as we noted in the introduction, define
\begin{align*}
S_{k \bl} =&\ g^{i \bj} \Omega_{i \bj k \bl}.
\end{align*}
Similarly define
\begin{align*}
P_{i \bj} =&\ g^{k \bl} \Omega_{i \bj k \bl}.
\end{align*}

\begin{lemma} \label{Scoordcalc} Given $(M^{2n}, g, J)$ a complex manifold one
has in complex
coordinates the formulas
\begin{align*}
S_{k \bl} =&\ g^{i \bj} \left( - g_{k \bl, i \bj} + g^{m \bn} g_{k \bn, i}
g_{\bl
m, \bj} \right)\\
P_{i \bj} =&\ g^{k \bl} \left( - g_{k \bl, i \bj} + g^{m \bn} g_{k \bn, i}
g_{\bl
m, \bj} \right).
\end{align*}
\begin{proof} These both follow from the general formula for the Chern curvature
given by
\begin{align*}
\Omega_{i \bj k \bl} =&\ - g_{k \bl, i \bj} + g^{m \bn} g_{k \bn, i} g_{\bl m,
\bj}.
\end{align*}
\end{proof}
\end{lemma}

\section{Existence and Regularity}
In this section we record some basic existence and regularity results for
solutions to (\ref{gflow}).  In particular we will see that this flow preserves
the pluriclosed condition, and is equivalent to (\ref{flow}) when the initial
condition is pluriclosed.  For notational convenience set
\begin{align} \label{phidef}
\Phi(\omega) = - \del \del^* \omega - \delb \delb^* \omega - \frac{\sqrt{-1}}{2}
\del \delb \log \det g.
\end{align}

\begin{prop} \label{ellipticity} Let $(M^{2n}, g, J)$ be a complex manifold with
pluriclosed metric.  Let $H = \{ g \in \Sym^2(T^*M) | g \mbox{ compatible with }
J, \del
\delb \gw =
0 \}$.  Then the operator
\begin{align*}
\omega \to \Phi(\gw)
\end{align*}
is a real quasi-linear second-order elliptic operator when restricted to $H$.
\begin{proof} Combining Lemmas \ref{operatorlemma1}, and
\ref{operatorlemma10} we see
\begin{align*}
\Phi(\omega)_{j \bk} =&\ - \frac{\sqrt{-1}}{2} g^{p \bq} \left[ g_{p \bk, \bq j}
+
g_{j \bq, p \bk} - g_{p \bq, j \bk}\right]\\
&\ + \frac{\sqrt{-1}}{2} g^{p \bm} g^{n \bq} \left[ g_{\bm n, j} \left( g_{p
\bk, \bq} - g_{p \bq, \bk} \right) + g_{n \bm, \bk} \left( g_{j \bq, p} - g_{p
\bq, j} \right) \right.\\
&\ \qquad \qquad \qquad \left. + g_{p \bq, j} g_{\bm n, \bk} \right]\\
=&\ - \frac{\sqrt{-1}}{2} \left[ g^{p \bq} \left(g_{p \bk, \bq j} +
g_{j \bq, p \bk} - g_{p \bq, j \bk} \right) \right.\\
&\ \left. \qquad \qquad + g^{p \bm} g^{n \bq} \left( g_{\bm n, j} g_{p \bk,
\bq} + g_{n \bm, \bk} g_{j \bq, p} - g_{p \bq, j} g_{\bm n, \bk} \right)\right].
\end{align*}
This coordinate formula shows that $\Phi$ is a real quasi-linear second-order
operator.  We now compute the symbol of the linearization to show that $\Phi$ is
elliptic when restricted to $H$.  Fix a point $x \in M$ and suppose the complex
coordinates are chosen such that
$g_{i \bj}(x) = \gd_{ij}$.  Moreover fix a family $g(s)$ of pluriclosed metrics
satisfying $\frac{\del}{\del s} g(s)_{|s = 0} = h$.  We can compute the
variation of $\Phi$ with respect to this family.  In particular we have
\begin{align} \label{linearization1}
\gs D \Phi(\omega)(h)_{j \bk} =&\ - \sum_{p = 1}^n h_{p \bk, \bp j} + h_{j \bp,
p
\bk} - h_{p \bp, j \bk}.
\end{align}
Note that the balancing condition for $g$ passes to the linearization.  In
particular for $i \neq j$ and $k \neq l$ we have
\begin{align*}
0 =&\ \left(\del \delb h( J \cdot, \cdot) \right)_{i j \bk \bl}\\
=&\ h_{i \bk, j \bl} + h_{j \bl, i \bk} - h_{i \bl, j \bk} - h_{j \bk, i \bl}.
\end{align*}
Suppose $j, k > 1$.  Applying the balancing condition we can write
\begin{align*}
h_{1 \bk, \bar{1}j} = h_{1 \bar{1}, j \bk} + h_{j \bk, 1 \bar{1}} - h_{j
\bar{1}, 1 \bk}.
\end{align*}
Therefore
\begin{align} \label{linearization2}
\gs D \Phi(\omega) (h)_{j \bk} =&\ - h_{j \bk, 1 \bar{1}} -
\sum_{p = 2}^n h_{p \bk, \bp j} + h_{j \bp, p \bk} - h_{p \bar{p}, j \bk}.
\end{align}
Now take the Fourier transform of the linear operator $D \Phi(\omega)$, and
further rotate coordinates such that the Fourier variable $\xi$ satisfies $\xi =
(1, 0 \dots, 0)$.  Now we can compute directly using (\ref{linearization1})
\begin{align*}
\left[\gs D \Phi (\omega) \right]^{\wedge}(h)(\xi)_{1 \bar{1}} =&\ - \brs{\xi}^2
h_{1 \bar{1}}.
\end{align*}
Likewise for $k > 1$ using (\ref{linearization1}) we compute
\begin{align*}
\left[\gs D \Phi (\omega) \right]^{\wedge}(h)(\xi)_{1 \bar{k}} =&\ -\brs{\xi}^2
h_{1 \bar{k}}.
\end{align*}
Finally for $j, k > 1$ we compute using (\ref{linearization2})
\begin{align*}
\left[\gs D \Phi (\omega) \right]^{\wedge}(h)(\xi)_{j \bar{k}} =&\ -\brs{\xi}^2
h_{j\bk}.
\end{align*}
Therefore $\Phi$ is elliptic, and the result follows.
\end{proof}
\end{prop}

Next we express equation (\ref{flow}) using the curvature and torsion of the
Chern connection.  Let $w$ denote the trace of the torsion.  In particular we
have in coordinates
\begin{gather*}
w_i = g^{j \bk} T_{i j \bk}.
\end{gather*}

\begin{prop} \label{torsionidentity} Let $(M^{2n}, g, J)$ be a complex manifold
with pluriclosed metric.  Then
\begin{align*}
\N w =&\ - \divg^{\N} T - Q^1.
\end{align*}
\begin{proof} Note that the pluriclosed condition implies
\begin{align*}
\del_{\bar{i}} T_{j k \bar{l}} = \del_{\bar{l}} T_{j k \bi}
\end{align*}
for any $i,j,k,l$.  We directly compute
\begin{align*}
\left( \N w \right)_{i \bj} =&\ g^{p \bq} \left( \N_{\bj} T_{i p \bq} \right)\\
=&\ g^{p \bq} \left( \del_{\bj} T_{i p \bq} - \gG_{\bj \bq}^{\br} T_{i p \br}
\right)\\
=&\ g^{p \bq} \left( \del_{\bq} T_{i p \bj} - \gG_{\bj \bq}^{\br} T_{i p \br}
\right)\\
=&\ g^{p \bq} \left( \N_{\bq} T_{i p \bj} + \left(\gG_{\bq \bj}^{\br} - \gG_{\bj
\bq}^{\br} \right) T_{i p \br} \right)\\
=&\ - \left(\divg^{\N} T\right)_{i \bj} - Q^1_{i \bj}.
\end{align*}
\end{proof}
\end{prop}

\begin{prop} \label{Sflow} Let $(M^{2n}, g, J)$ be a solution to (\ref{flow})
with pluriclosed initial condition.  Then
\begin{align*}
\dt g =&\ - S + Q^1.
\end{align*}
\begin{proof} Lemma \ref{operatorlemma1} implies that
\begin{align*}
\del \del^* \omega \left(J \cdot, \cdot \right) =&\ - \N w\\
\delb \delb^* \omega \left(J \cdot, \cdot \right) =&\ - \bar{\N w}.
\end{align*}
Thus composing the pluriclosed flow equation with $J$ yields
\begin{align*}
\dt g =&\ - \N w - \bar{\N w} - P.
\end{align*}
In general the Bianchi identity implies
\begin{align*}
P = S + \divg^{\N} T - \bar{\N w}.
\end{align*}
Applying the result of Proposition \ref{torsionidentity} yields
\begin{align*}
\divg^{\N} T =&\ - \N w - Q^1.
\end{align*}
Therefore
\begin{align*}
P = S - \N w - \bar{\N w} - Q^1
\end{align*}
and the result follows.
\end{proof}
\end{prop}

\begin{thm} Let $(M^{2n}, g, J)$ be a compact complex manifold with pluriclosed
metric $g$.  There exists a constant $c(n)$ depending only $n$ such that there
exists a unique solution $g(t)$ to (\ref{gflow})
for
\begin{align*}
t \in \left[0, \frac{c(n)}{\max \{ \brs{\Omega}_{C^0(g_0)}, \brs{\N
T}_{C^0(g_0)}, \brs{T}^2_{C^0(g_0)}
\}} \right].
\end{align*}
Moreover, there exist constants $C_m$ depending only on $m$ such that the
estimates
\begin{align*}
\brs{\N^m \Omega}_{C^0(g_t)}, \brs{\N^{m+1} T}_{C^0(g_t)} \leq \frac{C_m \max \{
\brs{\Omega}_{C^0(g_0)}, \brs{\N T}_{C^0(g_0)},
\brs{T}^2_{C^0(g_0)} \}}{t^{m / 2}}
\end{align*}
hold for all $t$ in the above interval.  Moreover, if $g(0)$ is pluriclosed the
metric $g(t)$ is
pluriclosed for
all $t$ and is a solution to (\ref{flow}).  If furthermore $g(0)$ is K\"ahler,
then $g(t)$ is K\"ahler for all
time and $g(t)$ solves K\"ahler-Ricci flow.
\begin{proof} The general
regularity theorem of \cite{ST} applies to give the existence
statements, and the claim that a K\"ahler initial condition results in a
solution to K\"ahler-Ricci flow.  For the convenience of the reader, we briefly
recall these arguments.  Lemma \ref{Scoordcalc} shows that
$S$ is a strictly elliptic operator, and so the short-time existence follows
from standard theory.  Furthermore, a judicious application of the Bianchi
identities shows that
\begin{gather} \label{curvatureev}
\begin{split}
\frac{\del}{\del t} \N^k \Omega =&\ \gD \N^k \Omega + \sum_{j = 0}^k \N^j T *
\N^{k + 1 - j}
\Omega + \sum_{j=0}^k \N^j
\Omega *
\N^{k-j} \Omega\\
&\ + \sum_{j=0}^{k} \sum_{l=0}^j \N^l T * \N^{j-l} T * \N^{k-j} \Omega.
\end{split}
\end{gather}
and
\begin{gather} \label{torsionev}
\begin{split}
\dt \N^k T =&\ \gD \N^k T + \sum_{j = 0}^{k+1} \N^j T * \N^{k+1 - j} T + \sum_{j
= 0}^{k} \N^j T * \N^{k - j} \Omega\\
&\ + \sum_{j = 0}^{k-1} \sum_{l = 0}^j \N^l T * \N^{j - l + 1} T * \N^{k - 1 -
j} T.
\end{split}
\end{gather}
One can apply standard estimates to these equations to derive the derivative
estimates.

We show that the pluriclosed condition is
preserved.  Suppose $g(0)$ is pluriclosed.  By Proposition \ref{ellipticity},
$\Phi$ is a strictly elliptic
operator, therefore (\ref{flow}) is a strictly parabolic equation and so
short-time existence and uniqueness of the solution to (\ref{flow}) with initial
condition $g(0)$ follows from standard results since $M$ is
compact.  We may directly compute to see that (\ref{flow}) preserves the
pluriclosed condition.  
Recall the equations $\del^2 = \delb^2 = 0$
and $\del \delb = -
\delb \del$.  Using this and the fact that the Chern form $\frac{\sqrt{-1}}{2}
\del \delb \log \det
g$ is closed, we directly compute
\begin{align*}
\dt \del \delb \omega =&\ - \del \delb \Phi(\omega)\\
=&\ \del \delb \left[ \left( \del \del^* \omega + \delb \delb^*
\omega \right) + \frac{\sqrt{-1}}{2} \del \delb \log \det g \right]\\
=&\ - \delb \del \del \del^* \omega + \del \delb \delb
\delb^* \omega\\
=&\ 0.
\end{align*}
Finally, using Proposition \ref{Sflow} we see that the solution to (\ref{flow})
in fact solves (\ref{gflow}).  Since solutions to (\ref{gflow}) are unique as
noted above, it follows that the solution to (\ref{gflow}) coincides with the
solution to (\ref{flow}), and therefore $g(t)$ is pluriclosed for all time.
\end{proof}
\end{thm}

\noindent A further consequence of the derivative estimates is a basic long-time
existence obstruction.
\begin{thm} \label{obstruction1} Let $(M^{2n}, g(t), J)$ be a solution to
(\ref{flow}).  Let $\tau$
denote the maximal existence time of the flow.  If $\tau < \infty$, then
\begin{align*}
\limsup_{t \to \tau} \max \left\{ \brs{\Omega}_{C^0(g_t)}, \brs{\N
T}_{C^0(g_t)}, \brs{T}^2_{C^0(g_t)} \right\} = \infty.
\end{align*}
\end{thm}

In the remainder of the section we compute evolution equations of basic integral
quantities and observe that in certain situations they function as monotonic
quantities along solutions to (\ref{flow}).
\begin{prop} Let $(M^{2n}, g(t), J)$ be a solution to (\ref{flow}) with
pluriclosed initial condition.  Then
the
volume of $g(t)$ satisfies
\begin{align*}
\dt \Vol(g(t)) =&\ 2 \int_M \brs{\del^* \omega}^2 - d.
\end{align*}
\begin{proof} We directly compute
\begin{align*}
\dt \Vol(g(t)) =&\ \dt \int_M dV_g\\
=&\ \int_M \tr_{\omega} \Phi(\omega) dV_g\\
=&\ \int_M \left< \omega, \del \del^* \omega + \delb \delb^* \omega +
\frac{\sqrt{-1}}{2} \del \delb \log \det g \right> dV_g\\
=&\ 2 \int_M \brs{\del^* \omega}^2 - d.
\end{align*}
\end{proof}
\end{prop}

Next we compute the evolution of the degree of a line bundle.  Recall the
definition of degree.
\begin{defn} Let $(M^{2n}, g, J)$ be a Hermitian manifold.  Let
\begin{gather} \label{ddef}
d = \deg(M) := \int_M \left< c_1(M), \omega \right> = \int_M \left( -
\frac{\sqrt{-1}}{2}
\del
\delb \log \det g \right) \wedge \omega^{n-1}.
\end{gather}
This is often called the degree of the surface.  More generally, given $\LL$ a
line bundle over $M$, define
\begin{gather} \label{degdef}
\deg(\LL) := \int_M c_1(\LL) \wedge \omega^{n-1}.
\end{gather}
\end{defn}

\begin{prop} \label{degreeev} Let $(M^4, g(t), J)$ be a solution to pluriclosed
flow on a complex
surface, and let $L$ be a line bundle over $M$.  Then
\begin{align*}
\dt \deg_{g_t}(L) =&\ - c_1(L)\cdot c_1(M).
\end{align*}
\begin{proof} We compute directly
\begin{align*}
\dt \deg_{g_t}(L) =&\ \dt \int_M c_1(L) \wedge \omega_t\\
=&\ \int_M c_1(L) \wedge \left( \del \del^* \omega + \delb
\delb^* \omega + \frac{\sqrt{-1}}{2} \del \delb \log \det g \right).
\end{align*}
Since $c_1(L)$ is closed,
\begin{align*}
\int_M c_1(L) \wedge \left(\del \del^* \omega + \delb \delb^* \omega
\right) = 0
\end{align*}
by Stokes Theorem.  Likewise
\begin{align*}
\int_M c_1(L) \wedge \frac{\sqrt{-1}}{2} \del \delb \log \det g =&\ - \int_M
c_1(L) \wedge c_1(M)\\
=&\ - c_1(L)\cdot c_1(M).
\end{align*}
\end{proof}
\end{prop}

\begin{prop} \label{areaev} Let $(M^{4}, g(t), J)$ be a solution to pluriclosed
flow on a
complex surface.  Fix $D$
a
divisor on $M$.  Then
\begin{align*}
\dt \int_D \omega =&\ c_1(D) \cdot c_1(M).
\end{align*}
\begin{proof} We directly compute, applying Stokes Theorem as in the previous
proposition,
\begin{align*}
\dt \int_D \omega =&\ \int_D \left( \del \del^* \omega +
\delb \delb^* \omega + \frac{\sqrt{-1}}{2} \del \delb \log \det g
\right)\\
=&\ \int_M \left( - c_1(D) \right) \wedge \left( \del \del^*
\omega + \delb \delb^* \omega + \frac{\sqrt{-1}}{2} \del \delb \log \det
g \right)\\
=&\ \int_M c_1(D) \wedge c_1(M)\\
=&\ c_1(D) \cdot c_1(M).
\end{align*}
\end{proof}
\end{prop}
Note that the time evolutions of the quantities in Propositions \ref{degreeev}
and \ref{areaev} are both topological in nature.  Indeed, the sign of this
topological number determines these quantities as monotonically increasing,
decreasing, or constant.  Once one has a sign for the degree, in certain cases
then the time evolution of volume is monotonic as well.  These same basic
observations are behind most of our classification results for static metrics in
the next section.

\section{Improved Regularity on Surfaces}

In this section we improve the basic regularity theory for (\ref{gflow}) in the
case of a complex surface.  
In the propositions below we give bounds for $\brs{T}^2$ and $\brs{\N T}$ in the
presence of a bound on $\brs{\Omega}$.  The key evolution equation appears in
Proposition \ref{Tnormev} at the end of this section.

\begin{prop} \label{torsionbound} Let $(M^{4}, g(t), J)$ be a solution to
(\ref{gflow}).  There exists a universal constant $c_0$ such
that if $\brs{\Omega} \leq C$ on $[0, \tau]$ then $\brs{T}^2 \leq \max \{
\brs{T}^2_{C^0(g_0)}, c_0 C \}$
on $[0, \tau]$.
\begin{proof} We will apply the maximum principle to the evolution equation for
$\brs{T}^2$.  According to Proposition \ref{Tnormev} we conclude that for a
solution to (\ref{gflow})
one has
\begin{align*}
\dt \brs{T}^2 =&\ \gD \brs{T}^2 - 2 \brs{\N T}^2 + \Omega * T^{*2} + \N
\brs{T}^2 * w - \frac{1}{2} \brs{T}^4.
\end{align*}
One can replace the term $\divg^{\N} T$ by Chern curvature terms using the
Bianchi identity.
Here the Laplacian is that of the Chern connection, however a calculation is
coordinates shows that
\begin{align*}
\gD f =&\ \gD_{LC} f - \left< w, \N f \right>.
\end{align*}
Applying the curvature bound it follows that
\begin{align*}
\dt \brs{T}^2 \leq&\ \gD_{LC} \brs{T}^2 + \N \brs{T}^2 * w + C \brs{T}^2 -
\frac{1}{2} \brs{T}^4.
\end{align*}
The result follows by the maximum principle.
\end{proof}
\end{prop}

\begin{prop} \label{torsionderbound} Let $(M^{4}, g(t), J)$ be a solution to
(\ref{gflow}).  Suppose the solution to (\ref{gflow}) exists on
$[0, T]$ and $\brs{\Omega} \leq K$ and $\brs{T}^2 \leq K$ on $[0, T]$.  There
exists a universal constant $C$ such that for any $\frac{1}{K} \leq t \leq T$
one has $\brs{\N T}(t) \leq C K$.
\begin{proof} We first recall some basic evolution inequalities.  In particular
from (\ref{curvatureev}) and (\ref{torsionev}) we conclude
\begin{align*}
\dt \brs{\N T}^2 \leq&\ \gD \brs{\N T}^2 + C \left( \brs{T} \brs{\N T} \brs{\N
\Omega} + \brs{\Omega} \brs{\N T}^2 + \brs{T}^2 \brs{\N T}^2 \right)\\
\dt \brs{\Omega}^2 \leq&\ \gD \brs{\Omega}^2 - 2 \brs{\N \Omega}^2 + C \left(
\brs{\Omega}^3 + \brs{T}^2 \brs{\Omega}^2 \right).
\end{align*}
We will derive an estimate for $\brs{\N T}$ at $t = \frac{1}{K}$ which implies
the general statement of the proposition.  Let
\begin{align*}
\Phi(x,t) =&\ t \left( \brs{\N T}^2 + \brs{\Omega}^2 \right) + A \brs{T}^2
\end{align*}
where $A$ is a constant to be chosen later.  Combining the evolution equations
above with Proposition \ref{Tnormev} and applying the assumed bound on curvature
and torsion we conclude
\begin{align*}
\dt \Phi =&\ \gD \Phi + C t \brs{T} \brs{\N T} \brs{\N \Omega} - 2 t \brs{\N
\Omega}^2 \\
&\ + \left( C t K \brs{\N T}^2 - 2 A \brs{\N T}^2 \right) + t K^3 + K^2. 
\end{align*}
Note that we may apply the Cauchy-Schwarz inequality to conclude
\begin{align*}
C t \brs{T} \brs{\N T} \brs{\N \Omega} - 2 t \brs{\N \Omega}^2 \leq C' t
\brs{T}^2 \brs{\N T}^2 \leq C' t K \brs{\N T}^2.
\end{align*}
As long as $t \leq \frac{1}{K}$ we can therefore choose $A$ large with respect
to universal constants so that
\begin{align*}
\dt \Phi \leq \gD \Phi + C K^2.
\end{align*}
Hence by the maximum principle we conclude $\Phi(\frac{1}{K}) \leq K$ and hence
\begin{align*}
\frac{1}{K} \brs{\N T}^2 \left(\frac{1}{K} \right) \leq C K
\end{align*}
and the result follows.
\end{proof}
\end{prop}

\begin{thm} Let $(M^4, g(t), J)$ be a solution to (\ref{gflow}).  Let $\tau$
denote the maximal existence time of the flow.  If $\tau < \infty$, then
\begin{align*}
\limsup_{t \to \tau} \brs{\Omega}_{C^0(g_t)} = \infty.
\end{align*}
\begin{proof} Suppose $\limsup_{t \to \tau} \brs{\Omega}_{C^0(g_t)} = C <
\infty$.  By Proposition \ref{torsionbound} the torsion is uniformly bounded up
to time $\tau$.  We may choose a small $\ge > 0$ so that Proposition
\ref{torsionderbound} applies for all $t \in [\ge, \tau]$ to yield a uniform
bound on $\brs{\N T}$ in this interval.  The result now follows from Theorem
\ref{obstruction1}.
\end{proof}
\end{thm}

In the remainder of this section we give a precise calculation of the evolution
of $\brs{T}^2$
along a solution to (\ref{gflow}).  Before we
begin we record a few algebraic identities which hold for quadratic expressions
in the torsion on a complex surface.  See (\ref{Qdefs}) for the definition of
$Q^2$.

\begin{lemma} \label{algebraidentities} Let $(M^4, g, J)$ be a Hermitian
surface.  Then
\begin{align*}
Q^1 =&\ \frac{1}{2} \brs{T}^2 g\\
\left< Q^2, Q^1 \right> =&\ \frac{1}{2} \brs{T}^4\\
\brs{Q^1}^2 =&\ \frac{1}{2} \brs{T}^4
\end{align*}
\begin{proof} Choose complex coordinates at a point so that $g$ is the identity
matrix.  On a surface there are two nonzero components of $T$ up to symmetry. 
In particular let
\begin{align*}
a = T_{1 2 \bar{1}}, \qquad b = T_{1 2 \bar{2}}.
\end{align*}
Then we can directly compute
\begin{align*}
Q^1 =&\ \left(
\begin{matrix}
a^2 + b^2 & 0\\
0 & a^2 + b^2
\end{matrix}
\right), \qquad
Q^2 = \left(
\begin{matrix}
2 a^2 & 2 ab\\
2 a b & 2 b^2
\end{matrix}
\right)
\end{align*}
All of the required identities follow immediately.
\end{proof}
\end{lemma}
\noindent Also, recall the Bianchi identities for the Chern connection.
\begin{lemma} \label{Bianchi} \emph{(Bianchi Identity)} Let $(M^{2n}, g, J)$ be
a Hermitian manifold and let $\N$ denote the Chern connection associated to $(g,
J)$.  For $X, Y, Z \in T_x(M)$
we have
\begin{align*}
\Sigma \{ \Omega(X, Y) Z \} =&\ \Sigma \{ T(T(X, Y), Z) + \N_X T(Y, Z) \}\\
\Sigma \{ \N_X \Omega(Y, Z) + \Omega(T(X, Y), Z) \} =&\ 0
\end{align*}
\end{lemma}

We start our calculation by using the general calculation for the evolution of
$T$:
\begin{lemma} \label{Tev} For a solution to $\dt g = - S + Q^1$ we have
\begin{align*}
\dt T_{i j \bk} =&\ \gD T_{i j \bk} + g^{m \bn} \left[T_{j i}^p \N_{\bn} T_{m p
\bk} + \N_{\bn} T_{m j}^p T_{i p
\bk} + T_{m j}^p \N_{\bn} T_{i p \bk} \right.\\
&\ \left. \qquad \quad + \N_{\bn} T_{i m}^p T_{j p \bk} + T_{i m}^p \N_{\bn}
T_{j p \bk} \right]\\
&\ + g^{m \bn} \left[\Omega_{\bn j m}^p T_{i p \bk} + \Omega_{\bn j \bk}^{\bp}
T_{i m \bp} - \Omega_{\bn i m}^p T_{j p \bk} \right.\\
&\ \qquad \quad \left. - \Omega_{\bn i \bk}^{\bp} T_{j m \bp} - \Omega_{p \bn m
\bk} T_{j i}^p \right] - T_{i j}^p \left(S_{p \bk} - Q^1_{p \bk} \right)\\
&\ + \N_i Q^1_{j \bk} - \N_j Q^1_{i \bk}.
\end{align*}
\begin{proof} This evolution equation is derived using the Bianchi identities
for the curvature of the Chern connection.  This calculation is found in
\cite{ST} Lemma 6.2.
\end{proof}
\end{lemma}

We would like to simplify this equation in the case of a solution to
(\ref{gflow}) on a surface.  In a series of lemmas below we simplify the inner
product of each term in the above lemma with
$T$.

\begin{lemma} \label{Tcalc1} Let $(M^{4}, g, J)$ be a Hermitian surface.  Then
\begin{align*}
g^{i \bm} g^{j \bn} g^{\bk p} \left( \N_i Q^1_{j \bk} - \N_j Q^1_{i \bk} \right)
T_{\bm \bn p} =&\ \left< \N \brs{T}^2, w \right>.
\end{align*}
\begin{proof} From Lemma \ref{algebraidentities} we know $Q^1 = \frac{1}{2}
\brs{T}^2 g$. 
Thus using metric compatibility of the connection we conclude
\begin{align*}
g^{i \bm} g^{j \bn} g^{\bk p} \left( \N_i Q^1_{j \bk} - \N_j Q^1_{i \bk} \right)
T_{\bm \bn p} =&\ \frac{1}{2} g^{i \bm} g^{j \bn} g^{\bk p} \left( g_{j \bk}
\N_i \brs{T}^2 - g_{i \bk} \N_j \brs{T}^2 \right) T_{\bm \bn p}\\
=&\ \left< \N \brs{T}^2, w \right>.
\end{align*}
\end{proof}
\end{lemma}

\begin{lemma} Let $(M^{4}, g, J)$ be a Hermitian surface.  Then
\begin{align*}
g^{i \bm} g^{j \bn} g^{\bk p} g^{r \bs} T_{j i}^q \left(\N_{\bs} T_{r q \bk} -
\Omega_{q \bs r \bk} \right) T_{\bm \bn p} =&\ \frac{1}{2} s \brs{T}^2.
\end{align*}
\begin{proof} We apply the Bianchi identity to conclude
\begin{align*}
\N_{\bs} T_{r q \bk} - \Omega_{q \bs r \bk} =&\ \Omega_{\bs r q \bk}.
\end{align*}
Plugging this in yields
\begin{align*}
g^{i \bm} g^{j \bn} g^{\bk p} g^{r \bs} T_{j i}^q \left(\N_{\bs} T_{r q \bk} -
\Omega_{q \bs r \bk} \right) T_{\bm \bn p} =&\ g^{i \bm} g^{j \bn} g^{\bk p}
g^{r \bs} T_{j i}^q \left( \Omega_{\bs r q \bk} \right) T_{\bm \bn p}\\
=&\ \left< S, Q^1 \right>\\
=&\ \frac{1}{2} \brs{T}^2 s.
\end{align*}
The result follows.
\end{proof}
\end{lemma}

\begin{lemma} Let $(M^{4}, g, J)$ be a Hermitian surface.  Then
\begin{align*}
g^{i \bm} g^{j \bn} g^{\bk p} g^{r \bs} \left( \left(\N_{\bs} T_{r j}^q +
\Omega_{\bs j r}^p \right) T_{i q \bk} - \left( \N_{\bs} T_{r i}^p + \Omega_{\bs
i r}^p \right) T_{j q \bk} \right)  T_{\bm \bn p} =&\ - s \brs{T}^2.
\end{align*}
\begin{proof} First we observe the symmetry between the pairs of terms and just
compute one of them.  Applying the Bianchi identity we conclude
\begin{align*}
\N_{\bs} T_{r j}^q + \Omega_{\bs j r}^q =&\ \Omega_{\bs r j}^q
\end{align*}
Plugging this in yields
\begin{align*}
g^{i \bm} g^{j \bn} g^{\bk p} g^{r \bs} \left(\N_{\bs} T_{r j}^q + \Omega_{\bs j
r}^p \right) T_{i q \bk}T_{\bm \bn p} =&\ g^{i \bm} g^{j \bn} g^{\bk p} g^{r
\bs} \Omega_{\bs r j}^q T_{i q \bk} T_{\bm \bn p}\\
=&\ - \left< S, Q^1 \right>\\
=&\ - \frac{1}{2} s \brs{T}^2.
\end{align*}
\end{proof}
\end{lemma}

\begin{lemma} \label{Tcalc4} Let $(M^{4}, g, J)$ be a Hermitian surface.  Then
\begin{align*}
g^{i \bm}& g^{j \bn} g^{\bk p} g^{r \bs} \left( T_{ir}^q \left( \N_{\bs} T_{j q
\bk} + \Omega_{\bs j \bk q} \right) - T_{j r}^q \left( \N_{\bs} T_{i q \bk} +
\Omega_{\bs i \bk q} \right) \right) T_{\bm \bn p}\\
=&\ \left< Q^2, S + \divg^{\N} T \right>.
\end{align*}
\begin{proof} We note that the two pairs of terms are the same using the
skew-symmetry of $T$.  We compute the first one.  Note
\begin{align*}
g^{i \bm} g^{j \bn} g^{\bk p} g^{r \bs} T_{ir}^q \left( \N_{\bs} T_{j q \bk} +
\Omega_{\bs j \bk q} \right) T_{\bm \bn p}=&\ g^{i \bm} g^{j \bn} g^{\bk p} g^{r
\bs} T_{ir}^q \left( \Omega_{j \bs q \bk} - \N_{\bs} T_{q j \bk} \right) T_{\bm
\bn p}\\
=&\ g^{i \bm} g^{j \bn} g^{\bk p} g^{r \bs} T_{i r}^q \Omega_{q \bs j \bk}
T_{\bm \bn p}.
\end{align*}
We want to reexpress this last term using the fact that $n = 2$.  We go through
the possibilities for the indices, assuming $g$ is the identity matrix at the
point we are computing at.  First suppose $i = 1$, then $m = 1$, $r = s = n =
j = 2$.  The resulting term is
\begin{align*}
g^{p \bk} T_{1 2}^q \Omega_{q \bar{2} 2 \bk} T_{\bar{1} \bar{2} p}.
\end{align*}
Likewise next assume $i = 2$, which implies $m = 1$, $r = s = n = j = 1$.  The
resulting term is
\begin{align*}
g^{p \bk} T_{2 1}^q \Omega_{q \bar{1} 1 \bk} T_{\bar{2} \bar{1} p}.
\end{align*}
So, let
\begin{align*}
A_{i \bj} = g^{k \bl} \Omega_{k \bj i \bl}.
\end{align*}
It is clear from the above calculations that
\begin{align*}
g^{i \bm} g^{j \bn} g^{\bk p} g^{r \bs} T_{i r}^q \Omega_{q \bs j \bk} T_{\bm
\bn p} =&\ \frac{1}{2} \left< Q^2, A \right>.
\end{align*}
Finally we want to apply the Bianchi identity once more to simplify the tensor
$A$.  In particular we see
\begin{align*}
A_{i \bj} =&\ g^{k \bl} \Omega_{k \bj i \bl}\\
=&\ g^{k \bl} \Omega_{\bj k \bl i}\\
=&\ g^{k \bl} \left( \Omega_{\bl k \bj i} + \N_{k} T_{\bl \bj i} \right)\\
=&\ S_{i \bj} + \divg^{\N} T_{\bj i}.
\end{align*}
The result now follows.
\end{proof}
\end{lemma}

\begin{prop} \label{Tnormev} Let $(M^4, g(t), J)$ be a complex surface with
$g(t)$ a solution to (\ref{gflow}).  Then
\begin{align*}
\dt \brs{T}^2 =&\ \gD \brs{T}^2 - 2 \brs{\N T}^2 + 2 \left<\N \brs{T}^2, w
\right> + \left<
Q^2, S + 2 \divg^{\N} T \right> - \frac{1}{2} \brs{T}^4
\end{align*}
\begin{proof} We start with the basic calculation
\begin{align*}
\dt \brs{T}^2 =&\ 2 \left< \dt T, T \right> + \left< 2 Q^1 + Q^2, S - Q^1
\right>.
\end{align*}
Next we plug the results of Lemmas \ref{Tcalc1} - \ref{Tcalc4} into Lemma
\ref{Tev} and apply Lemma \ref{algebraidentities} to conclude
\begin{align*}
\dt \brs{T}^2 =&\ 2 \left< \gD T, T \right> + 2 \left< \N \brs{T}^2, w \right>
- s \brs{T}^2\\
&\ + 2 \left< Q^2, S + \divg^{\N}  T \right> - 2 \left<Q^2, S - Q^1 \right> +
\left< 2 Q^1 + Q^2, S - Q^1 \right>\\
=&\ \gD \brs{T}^2 - 2 \brs{\N T}^2 + 2 \left<\N \brs{T}^2, w \right> + \left<
Q^2, S + 2 \divg^{\N} T \right>\\
&\ + \left< Q^2, Q^1 \right> - 2 \brs{Q^1}^2\\
=&\ \gD \brs{T}^2 - 2 \brs{\N T}^2 + 2 \left<\N \brs{T}^2, w \right> + \left<
Q^2, S + 2 \divg^{\N} T \right> - \frac{1}{2} \brs{T}^4
\end{align*}
as required.
\end{proof}
\end{prop}

\section{Static Metrics on Surfaces} \label{static}
In this section we derive certain identities satisfied by \static metrics.  As
applications of these identities we will classify static metrics on K3 surfaces,
two-dimensional complex tori, surfaces of general type, and nonprimary Hopf
surfaces.

\begin{defn} Let $(M^{2n}, g, J)$ be a complex manifold with pluriclosed metric.
 We say that $g$ is \emph{\static}if
\begin{align*}
\Phi(\omega) = \gl \omega
\end{align*}
for some constant $\gl$, and
\begin{align*}
\Vol(g) = 1.
\end{align*}
\end{defn}
\noindent We have ruled out the scaling ambiguity by fixing the volume to be
$1$.

\begin{prop} \label{staticidentity1} Let $(M^{4}, g, J)$ be a complex surface
with a \static metric.  Then
\begin{align*}
d - 2 \gl = 2 \int_M \brs{\del^* \omega}^2.
\end{align*}
\begin{proof} We begin by taking the wedge product of the \static equation with
$\omega$
and integrating.  Since $\int_M \omega \wedge \omega = 2$ this yields
\begin{align*}
- 2 \gl =&\ \int_M \left( \del \del^* \omega + \delb \delb^*
\omega + \cf  \right) \wedge \omega\\
=&\ \int_M \left< \left( \del \del^* \omega + \delb \delb^*
\omega + \cf  \right) , \omega \right>\\
=&\ \int_M \brs{\del^* \omega}^2 + \int_M \brs{\delb^*
\omega}^2 - d.
\end{align*}
Since $\brs{\del^* \omega}^2 = \brs{\delb^* \omega}^2$ the result follows.
\end{proof}
\end{prop}

\begin{prop} \label{staticidentity2} Let $(M^{4}, g, J)$ be a complex surface
with \static metric, and let $L$ be a line bundle on $M$.  Then
\begin{align*}
c_1(M) \cdot c_1(L) = \gl \deg L.
\end{align*}
\begin{proof} Let $\Omega$ be a (closed) form representing $c_1(L)$.  Take the
wedge product of the \static equation with $\Omega$ and integrate.  This yields
\begin{align*}
\int_M \left( - \del \del^* \omega - \delb \delb^* \omega - \frac{\sqrt{-1}}{2}
\del \delb \log \det g \right) \wedge \Omega =&\
\gl \int_M \omega \wedge \Omega.
\end{align*}
By definition the integral on the right hand side is the degree of $L$.  As for
the left hand side, since $\Omega$ is closed
\begin{align*}
\int_M \left( \del \del^* \omega + \delb \delb^* \omega \right)
\wedge \Omega = 0
\end{align*}
by Stokes Theorem.  The remaining term is $c_1(M) \cdot c_1(L)$.
\end{proof}
\end{prop}
There is a further quadratic identity we have for \static metrics.  First we
require a certain reverse H\"older type inequality (\cite{Buchdahl} Lemma 4). 
We include the proof for the reader's
convenience.
\begin{lemma} (\cite{Buchdahl} Lemma 4) \label{auxlemma}  If $\psi \in
\Lambda^{1,1}_{\mathbb R}$
satisfies $\del \delb \psi = 0$, then
\begin{align*}
\left( \int_M \omega \wedge \psi \right)^2 \geq \left( \int_M \omega^2 \right)
\left( \int_M \psi^2 \right),
\end{align*}
and moreover equality holds if and only if $\psi = c \omega + i \del \delb g$
for some constant $c$ and some $g \in C^{\infty}(M)$.
\begin{proof} Let $c := \frac{\int_M \omega \wedge \psi}{\int_M \omega^2}$. 
Since the function solving $f dV = \omega \wedge \psi - c \omega$ is
$L^2$-orthogonal to the constants, which are the only kernel of $f \to * i \del
\delb (\omega f)$, it follows that there exists a function $\rho$ such that
\begin{align*}
\omega \wedge \left( \psi - c \omega - i \del \delb \rho \right) = 0.
\end{align*}
It follows that the form $\psi - c \omega - i \del \delb \rho$ is antiselfdual,
hence
\begin{align*}
0 \leq \brs{\brs{\psi - c \omega - i \del \delb \rho}}^2 =&\ - \int_M \left(
\psi - c \omega - i \del \delb g \right)^{\wedge 2}\\
=&\ - \int_M \psi^{\wedge 2} + \frac{\left( \int_M \psi \wedge \omega
\right)^2}{\int_M \omega^2}.
\end{align*}
The proposition follows from this inequality.
\end{proof}
\end{lemma}

\begin{prop} \label{staticidentity3} Let $(M^4, J, g)$ be a compact complex
surface with \static pluriclosed metric.  Then
\begin{align*}
c_1^2 - 2 \gl d + \frac{1}{2} d^2 \geq&\ 0,\\
- c_1^2 + \frac{1}{2} d^2 \geq 0
\end{align*}
with equality in either case if and only if either $g$ is K\"ahler-Einstein or
$c_1(M) = 0$.
\begin{proof} We directly compute
\begin{align*}
2 \gl^2 =&\ \int_M \gl \omega \wedge \gl \omega\\
=&\ \int_M \left( \del \del^* \omega + \delb \delb^* \omega +
\frac{\sqrt{-1}}{2} \del \delb \log \det g \right)^{\wedge 2}\\
=&\ c_1^2(M) + \int_M \left( \del \del^* \omega + \delb \delb^* \omega
\right)^{\wedge 2}
\end{align*}
The last line follows by Stokes Theorem since $\frac{\sqrt{-1}}{2} \del \delb
\log \det g$ is
closed.  Now we apply Lemma \ref{auxlemma} to conclude
\begin{align*}
\int_M \left( \del \del^* \omega + \delb \delb^* \omega
\right)^{\wedge 2} \leq&\ \frac{1}{2} \left( \int_M \omega \wedge \left( \del
\del^* \omega + \delb \delb^* \omega \right) \right)^2\\
=&\ 2 \left( \int_M \brs{\del^* \omega}^2 \right)^2\\
=&\ \frac{1}{2} \left(2 \gl - d \right)^2\\
=&\ 2 \gl^2 - 2 \gl d + \frac{1}{2} d^2.
\end{align*}
The second to last line follows from Proposition \ref{staticidentity1}. 
Plugging this into the above yields the first inequality, and the second follows
using $\gl d = c_1(M)^2$, which follows from Proposition
\ref{staticidentity2}.
To characterize the equality case we note by Lemma \ref{auxlemma} that it occurs
if and only if
\begin{align*}
\del \del^* \omega + \delb \delb^* \omega = c \omega
+ \del \delb g.
\end{align*}
Plugging this into the \static equation yields
\begin{align*}
\frac{\sqrt{-1}}{2} \del \delb \log \det g + \del \delb g =&\ - \left(\gl + c
\right) \omega.
\end{align*}
The left hand side is closed, so if $\gl + c \neq 0$ then $\omega$ is K\"ahler,
hence K\"ahler-Einstein.  Otherwise we have that $\frac{\sqrt{-1}}{2} \del \delb
\log \det g$ is
represented by an exact form, therefore $c_1(M) = 0$.
\end{proof}
\end{prop}

\begin{prop} Consider a K3 surface or torus $(M^4, g, J)$ with \static
metric $g$.  Then $g$ is in fact K\"ahler and Ricci-flat.
\begin{proof}  Note that $c_1(M) = 0$.  Since $M$ has a K\"ahler Ricci-flat
metric $g_0$, it follows that
\begin{align*}
d(M, g) =&\ \int_M - \frac{\sqrt{-1}}{2} \del \delb \log \det g \wedge \omega\\
=&\ \int_M - \frac{\sqrt{-1}}{2} \del \delb \left( \log \frac{\det g}{\det g_0}
+
\log \det g_0 \right) \wedge \omega\\
=&\ 0.
\end{align*}
The last follows by integrating by parts, using that $\del \delb \omega = 0$.
Also, the K\"ahler form of $g_0$, $\omega_0$, is closed.  Taking the
wedge product of the \static equation with $\omega_0$ and integrating yields
\begin{align*}
\gl \int_M \omega \wedge \omega_0 =&\ \int_M \left( \del \del^* \omega +
\delb \delb^* \omega + \frac{\sqrt{-1}}{2} \del \delb \log \det g
\right) \wedge \omega_0\\
=&\ 0.
\end{align*}
However, since $\omega$ and $\omega_0$ are both positive $(1,1)$ forms one has
$\int_M \omega \wedge \omega_0 > 0$, therefore $\gl =
0$.  Therefore by Proposition \ref{staticidentity1} we conclude
\begin{align*}
\int_M \brs{\del^* \omega}^2 = d - 2 \gl = 0,
\end{align*}
therefore $g$ is K\"ahler and hence K\"ahler-Einstein.  Since $c_1(M) = 0$ this
means $g$ is Ricci-flat.
\end{proof}
\end{prop}

\begin{prop} \label{generaltype} Let $(M^{4}, g, J)$ be a surface of general
type with static metric.  Then $g$ is K\"ahler-Einstein.
\begin{proof} Since $(M^4, J)$ is of general type we have $c_1^2 > 0$, and also
one has $d < 0$ since the canonical bundle is ample.  It follows from
Proposition \ref{staticidentity2} that $\gl < 0$.  By Proposition
\ref{staticidentity3} we
conclude $d^2 \geq 2 c_1^2$, therefore $d \leq - \sqrt{2 c_1^2}$.  Returning to
Proposition \ref{staticidentity2} we conclude
\begin{align*}
c_1^2 =&\ \gl d \geq - \gl \sqrt{2 c_1^2}.
\end{align*}
Therefore $\gl \geq - \sqrt{\frac{c_1^2}{2}}$.  It follows that
\begin{align*}
d - 2 \gl \leq&\ - \sqrt{2 c_1^2} + 2 \sqrt{\frac{c_1^2}{2}} = 0.
\end{align*}
Now it follows from Proposition \ref{staticidentity1} that $\del^* \omega = 0$
therefore $\del \omega = 0$ and the metric is K\"ahler, hence K\"ahler-Einstein.
\end{proof}
\end{prop}

\begin{prop} \label{curves} Let $(M^4, g, J)$ be a complex surface with \static
metric and
suppose $\Sigma \subset M$ is a holomorphic curve such that $[\Sigma] = 0 \in
H^2(M, \mathbb R)$.  Then $\gl = 0$.
\begin{proof} We simply integrate the \static equation along the curve $\Sigma$
to yield
\begin{align*}
\gl \int_{\Sigma} \omega =&\ \int_{\Sigma} \left(- \del \del^* \omega - \delb
\delb^* \omega - \frac{\sqrt{-1}}{2} \del \delb \log \det g \right)\\
=&\ \int_M c_1(\Sigma) \wedge \left( - \del \del^* \omega - \delb \delb^* \omega
-
\frac{\sqrt{-1}}{2} \del \delb \log \det g \right)\\
=&\ c_1(\Sigma) \cdot c_1(M).
\end{align*}
The last line follows applying Stokes Theorem since $c_1(\Sigma)$ is closed. 
However, since $\Sigma$ is null-homologous, $c_1(\Sigma) = 0$, therefore the
right hand side is zero.  Since $\int_{\Sigma} \omega > 0$, therefore $\gl = 0$.
\end{proof}
\end{prop}

\begin{prop} Let $(M^{4}, J)$ be a Hopf surface blown up at $p > 0$ points. 
Then $M$ admits no static metrics.
\begin{proof} We start by observing that any Hopf surface admits a holomorphic
curve which is homologous to zero.  A general Hopf surface is defined by taking
the quotient of $\mathbb C^2 \setminus \{0\}$ by the group $G$ generated by the
map
\begin{align*}
\left(z_1, z_2 \right) \to \left(\ga_1 z_1, \ga_2 z_2 \right)
\end{align*}
where $0 < \brs{\ga_1} \leq \brs{\ga_2} < 1$.  It follows from \cite{BPV}
Proposition 18.2 that the space $H_{\ga} = \left( \mathbb C^2 \setminus \{0\}
\right) / G$ is either an elliptic fibration over $\mathbb P^1$ or contains
exactly two irreducible curves.  Thus any case contains at least one holomorphic
curve, and since all $H_{\ga}$ are homeomorphic to $S^3 \times S^1$ this curve
is necessarily null-homologous since $H_2(S^3 \times S^1, \mathbb R) = 0$.  It
is clear that this curve still exists, and is still null-homologous, if we
blow-up the Hopf surface at finitely many points.  It follows from Proposition
\ref{curves} that $\gl = 0$ for our \static metric.  By Proposition
\ref{staticidentity2} we yield $c_1^2(M) = \gl d = 0$.  However, if we have
blown up at $p$ points we have $c_1^2(M) = -p$, a contradiction.  Thus the
result follows.
\end{proof}
\end{prop}

We now come to a key observation on \static metrics.  In particular, we find
that manifolds admitting static metrics with nonzero constant admit
Hermitian-symplectic
structures, which were defined in the introduction.

\begin{prop} \label{symplecticexistence} Let $(M^{2n}, g, J)$ be a compact
complex manifold with \static
metric.  If $\gl \neq 0$, then $M$ is a Hermitian-symplectic manifold and
specifically
$\omega$ is the $(1,1)$-part of a Hermitian-symplectic form $\til{\omega}$. 
Furthermore, one has
\begin{align*}
\int_M \til{\omega}^{\wedge n} > 0.
\end{align*}
\begin{proof} Consider the real two-form
\begin{align*}
\til{\omega} =&\ \gw - \frac{1}{\gl} \left(\delb \del^* \omega + \del \delb^*
\omega
\right).
\end{align*}
We compute using the \static equation
\begin{align*}
d \til{\gw} =&\ d \left( \omega - \frac{1}{\gl} \left( \delb \del^* \omega +
\del \delb^* \omega \right) \right)\\
=&\ - \frac{1}{\gl} d \left( d \del^* \omega + d \delb^* \omega +
\frac{\sqrt{-1}}{2} \del \delb \log \det g \right)\\
=&\ 0.
\end{align*}
Thus $\til{\omega}$ is a Hermitian-symplectic form.  The claim that $\int_M
\til{\omega}^{\wedge n} > 0$ is a standard fact whose proof we outline below,
following calculations of \cite{Delanoe}.  In particular, set
\begin{align*}
\til{\omega} =&\ \phi + H + \bar{\phi}
\end{align*}
where $H  = \omega$ is the $(1,1)$ part of $\til{\omega}$ and $\phi$ is the
$(2,0)$ part of $\til{\omega}$, with $\bar{\phi}$ then being the $(0,2)$ part. 
Let $\{\theta^j\}$ denote a basis for $(1,0)$ forms at a point $x \in M$ which
diagonalizes the metric $\omega$.  We note
\begin{align*}
\til{\omega}^{\wedge n} =&\ \sum_{k = 0}^{\left[\frac{n}{2}\right]}
\frac{n!}{(k!)^2 (n - 2k)!} i^r H_{j_1 \bar{j}_1} \dots H_{j_r \bar{j}_r}
\cdot\\
&\ \qquad 
\theta^{j_1} \wedge \bar{\theta}^{j_1} \wedge \dots \wedge \theta^{j_r} \wedge
\bar{\theta}^{j_r} \wedge \phi^{\wedge k} \wedge \bar{\phi}^{\wedge k}.
\end{align*}
Note that the term $k = 0$ in the above summand is \emph{strictly} positive
since all of the coefficients $H_{j \bar{j}}$ are strictly positive and the
resulting form is then a positive multiple of the volume form of $\omega$.  It
remains to show that the terms $k > 0$ are nonnegative.  Fix a particular term
in the summand, and assume (without loss of generality by relabelling) that
$\{j_i = 2k + i\}$.  One can directly compute (\cite{Delanoe} pg. 845) that
\begin{align*}
\phi^{\wedge k} \wedge \bar{\phi}^{\wedge k} = i^{2k} f \theta^{1} \wedge
\theta^{\bar{1}} \wedge \dots \wedge \theta^{2k} \wedge \bar{\theta}^{2k}.
\end{align*}
where $f \geq 0$.  Thus the terms in the above summand with $k > 0$ are
nonnegative multiples of the volume form of $\omega$, and the result follows.
\end{proof}
\end{prop}

\section{Class VII Surfaces}
A minimal compact complex surface $S$ is of Kodaira's Class $\seven$ if $b_1(S)
= 1$.  The standard examples of surfaces of class $\seven$ with $b_2(S) = 0$ are
the Hopf surfaces, defined in section 4.  As we note in the following
example, the standard Hopf surface admits a static metric with $\gl = 0$.

\begin{ex} Consider the quotient of $\mathbb C^2 \setminus \{0\}$ by the group
$\Gamma$ generated by the map
\begin{gather*}
(z_1,z_2) \to \left( \frac{1}{2} z_1, \frac{1}{2} z_2 \right).
\end{gather*}
The standard complex structure on $\mathbb C^2$ descends to the quotient
manifold, which is homeomorphic to $S^3 \times S^1$.  Furthermore, setting $\rho
= \sqrt{z_1\bar{z}_1 + z_2 \bar{z}_2}$, the K\"ahler form
\begin{gather*}
\omega = \frac{1}{\rho^2} \del \delb \rho^2
\end{gather*}
is $\Gamma$-invariant and compatible with the standard complex structure, and
therefore defines a metric on the quotient.  This metric is pluriclosed,
although it should be noted that this is not true for higher dimensional Hopf
surfaces.  Call the associated metric $g$.  One can directly compute that
\begin{gather*}
S_g = g
\end{gather*}
and
\begin{gather*}
Q^1_g = g
\end{gather*}
so that
\begin{gather*}
S_g - Q^1_g = 0.
\end{gather*}
Thus $g$ is a static metric with $\gl = 0$.  This constant is of course
correctly predicted by Proposition \ref{curves}.
\end{ex}

It may be the case that every class $\seven$ surface with $b_2(S) = 0$ admits a
static metric.  The situation when $b_2(S) > 0$ is completely different though,
with none admitting static metrics. We say that a class $\seven$ surface is in
class $\seven^+$ if $b_2(S) > 0$.  We recall a basic
structural theorem for class $\seven^+$ surfaces.
\begin{thm} \label{structure} (\cite{Dloussky1} Theorem 1.8) Let $S$ be a (not
necessarily
minimal) compact complex surface such that $b_1(S) = 1$, and $b_2(S) = n > 0$. 
Then there exist $n$ exceptional line bundles $L_j, j = 0, \dots, n-1$, unique
up to torsion by a flat line bundle $F \in H^1(S, \mathbb C^*)$ such that
\begin{itemize}
\item{ $E_j = c_1(L_j), 0 \leq j \leq n - 1$ is a $\mathbb Z$-basis of $H^2(S,
\mathbb Z)$.}
\item{ $K_S L_j = -1$ and $L_i L_j = - \gd_{j k}$}
\item{ $K_S = L_0 + \dots + L_{n - 1} \in H^2(M, \mathbb Z)$}
\end{itemize}
\begin{proof} In fact we have only stated the part of (\cite{Dloussky1} Theorem
1.8) which we will need.  We just sketch a couple of the ideas.  One applies
Donaldson's
Theorem \cite{Donaldson} to diagonalize the intersection form on the
torsion-free projection of $H^2(M, \mathbb
Z)$.  Since the Kodaira dimension of a class $\seven^+$ surface is
$-\infty$, we have
that the geometric genus $p_g(M) = 0$.  By \cite{BPV} Theorem 2.7 (iii) we know
that $b^+(M) = 2 p_g(M) = 0$.  Therefore $b^+(M) = 0$ and the intersection form
of $M$ is negative definite.  Since $p_g = h^2(S, \mathcal O_S) = 0$, each
element of the basis is
realized as $c_1(L)$, and so the first claim follows.  The second and third and
fourth follow from applications of Riemann-Roch.  For the proof of the final
claim see \cite{Dloussky1}.
\end{proof}
\end{thm}

\begin{prop} \label{class7prop} Class $\seven^+$ surfaces admit no \static
metrics.
\begin{proof} Let $n = b_2(M) > 0$.  It follows from Theorem \ref{structure}
that
\begin{align*}
c_1(M)^2 = K_S \cdot K_S = -n.
\end{align*}
Therefore it follows from Proposition \ref{staticidentity2} that
\begin{align*}
-n = \gl d.
\end{align*}
Therefore $\gl \neq 0$.  Proposition \ref{symplecticexistence} implies that $M$
carries the structure of a symplectic manifold given by $\til{\omega}$. 
Moreover, $\int_M \til{\omega}^{\wedge 2} > 0$, so that $\til{\omega}$ is a
cohomology class with positive self-intersection.  However, as we remarked in
the proof of Theorem \ref{structure}, the intersection form on $M$ is negative
definite.  Thus we have arrived at a contradiction.
\end{proof}
\end{prop}

Proposition \ref{class7prop} is an important first step in applying solutions to
(\ref{gflow})
to studying the topology of class $\seven^+$ surfaces.  Recall that surfaces of
class $\seven$ with $b_2 = 0$ are classified.  In particular they are
biholomorphic to either a Hopf surface or an Inoue surface (a free quotient of
$\mathbb C \times \mathbb H$ by a properly discontinuous affine action).  This
result is typically called Bogomolov's Theorem \cite{Bo1}, \cite{Bo2},
although the first complete proofs appear to have been given independently in
\cite{LYZ} and \cite{Teleman}.  As for surfaces of class $\seven^+$, (i.e. $b_2
> 0$), in deep recent work Teleman \cite{Tel2}, \cite{Tel3} has confirmed the
conjectural list for $b_2 = 1, 2$ using Donaldson theory.  The dedicated effort
of many authors, culminating in the theorem of Dloussky-Oeljeklaus-Toma
\cite{DOT}, has reduced the problem of classifying class $\seven^+$ surfaces to
finding $b_2$ rational curves in a given minimal surface of class $\seven^+$,
and this is what Teleman shows for $b_2 = 1, 2$.  What Theorem \ref{class7} says
is that equation (\ref{gflow}) must encounter some nontrivial singularities on a
class $\seven^+$ surface.  In principle, these singularities should be closely
related to curves on the surface.  Indeed, such behavior is seen for the
K\"ahler-Ricci flow, where the flow exhibits blowdowns along complex
subvarieties \cite{TiZha}.  If sufficiently many curves can be shown to exist as
singularities of (\ref{gflow}), one could finish the classification
of class $\seven$ surfaces.  Much work remains to see to what extent
singularities of (\ref{gflow}) can be analyzed.

\bibliographystyle{hamsplain}

\end{document}